\documentclass[11pt]{article}
\usepackage[small]{titlesec}
\usepackage[top = 0.66in,textwidth = 6.5in, textheight=9.1in]{geometry}

\usepackage{amsmath}
\usepackage{graphicx}
\usepackage{latexsym}
\usepackage{color}
\usepackage{amssymb}
\usepackage{accents}
\usepackage{tabularx}
\usepackage{fancyhdr}
\usepackage{verbatim}
\usepackage{multirow}
\usepackage{framed}
\usepackage[square,sort,comma,numbers]{natbib}
\usepackage{url}
\usepackage{float, subfig}
\usepackage{enumitem}
\usepackage{mathtools}
\usepackage{mathrsfs}
\usepackage{amsfonts}
\usepackage{listings}
\usepackage{amsthm}
\usepackage{grffile}
\usepackage{sidecap}
\usepackage{pbox}
\usepackage[english,ruled,vlined,linesnumbered]{algorithm2e}
\usepackage{longtable}
\usepackage[noend]{algpseudocode}
\usepackage{scalefnt}
\usepackage{nomencl}
\usepackage{etoolbox}
\usepackage{diagbox}
\usepackage{titling}

\def\qed{\hfill{\(\vcenter{\hrule height1pt \hbox{\vrule width1pt height5pt
     \kern5pt \vrule width1pt} \hrule height1pt}\)} \medskip}

\newcommand{\bE}{{\mathbb{E}}}

\newcommand{\cX}{{\cal X}}

\newcommand{\noi}{\noindent}

\newcommand{\ms}{\medskip}

\newcommand{\st}{\mbox{s.t.}}

\makenomenclature
\renewcommand\nomgroup[1]{%
    \item[\bfseries
    \ifstrequal{#1}{A}{Indices and index sets}{%
    \ifstrequal{#1}{B}{Parameters}{%
    \ifstrequal{#1}{C}{Decision variables}{%
    }}}%
]\vspace{0.1in}}

\allowdisplaybreaks

\usepackage{authblk}
\title{\textbf{A Review of Sequential Decision Making via Simulation}}
\author[b]{Zhuo Zhang}
\author[a]{Dan Wang}
\author[b]{Haoxiang Yang}
\author[a,*]{Shubin Si}

\affil[a]{\textit{Department of Industrial Engineering, School of Mechanical Engineering, Northwestern Polytechnical University, Xi’an 710072, China}}
\affil[b]{\textit{School of Data Science, the Chinese University of Hong Kong, Shenzhen, Shenzhen, 518172, China}}
\affil[*]{\textit{Corresponding author: sisb@nwpu.edu.cn (S. Si).}}

\begin{document}
\baselineskip0.25in
\maketitle

\begin{abstract}
    Optimization via simulation has been well established to find optimal solutions and designs in complex systems. However, it still faces modeling and computational challenges when extended to the multi-stage setting. This survey reviews the models and methodologies of single-stage optimization via simulation and multi-stage stochastic programming. These are necessary theoretical components to push forward the development of sequential decision making via simulation. We identify the key challenge of sequential decision making via simulation as the appropriate modeling of the stage-wise value function, for which we survey the state-of-the-art meta-models and their potential solution algorithms.
\end{abstract}

\noi \textbf{\emph{Keywords}:} optimization under uncertainty; sequential decision process; multi-stage stochastic program; optimization via simulation.

\section{Introduction} \label{sec:intro}
    Optimization via simulation refers to using simulation techniques and methodologies to optimize complex systems or processes, especially under the consideration of uncertainty. For many applications, it is difficult to directly obtain the analytical form of the performance measure function in terms of decision variables, i.e., the system configuration. The decision maker makes use of computer-based or physical simulation models to analyze and evaluate such performance measures under different scenarios and then identifies the optimal or near-optimal decision variables that result in the best performance according to predefined criteria, such as expected cost, equity, reliability, and other utilities. There is usually a realistic limit for the round of simulation due to a physical budget, for example, the test of engine design~\cite{koegeler2000using}, or a computational limit in time and resources~\cite{chen2011stochastic}. The objective of the simulation optimization often lies within the following three categories: 1) to achieve the smallest regret compared to the optimal policy in the long run, under the assumption that the whole system dynamics are known; 2) to maximize the probability of choosing the optimal system configuration under a fixed simulation budget; and 3) to minimize the number of simulation runs to select the optimal solution with a fixed probabilistic guarantee.   
    
    Optimization via simulation has gained attention in both simulation and optimization societies because of its broad applicability in many real-world problems; for example, chemical engineering processes where the fundamental chemical and physical relationships can be well approximated by simulation processes~\cite{biegler1997systematic,bhosekar2018advances,biegler2014recent}, transportation and supply chains where simulation models characterize the flow dynamics and agent behavior~\cite{oliveira2016perspectives, romero2012simulation}, drug design where molecular dynamics simulation assists optimizing the drug structure and the combination of small compounds~\cite{alonso2006combining, yu2017computer}, manufacturing processes where complex material flow relationships are captured by simulation models to help determine optimal machine schedules~\cite{ramasesh1990dynamic, van1992job, liu2018survey, liu2011production}, and pandemic control where it is natural to use a simulation model in the place of the differential equations representing pandemic transmission dynamics~\cite{nsoesie2013simulation, paleshi2017simulation, yang2021design}. A common feature of the applications mentioned above of optimization via simulation methods is that since it is hard to find a closed-form representation of many processes in science, engineering, and human society, researchers resort to characterizing a simulation system to approximate the reality and hope to make the optimal decision based on the feedback from the simulation system. 

    We can formulate a general model of optimization via simulation as follows:
    \begin{align}
        \min_{x \in \cX} \quad & \bE_\xi[f(x,\xi)], \label{prob:stoch_prob}
    \end{align}
    where we aim to find the optimal solution \(x\) within a feasible region \(\cX\). Even with the same solution \(x\), the uncertainty from the environment \(\xi\) will affect the value of the objective function \(f\), and such uncertainty may or may not depend on the choice of \(x\). Since we do not know the realization of the uncertainty \(\xi\) when we make the decision, we aim to optimize the expected value. The expected value operator can be substituted by other utility functions such as conditional value at risk (CVaR) in risk-averse optimization~\cite{makarova2021risk} and maximization in robust optimization~\cite{bogunovic2018adversarially}. 

    When the closed form of \(f\) is not known and evaluating such a function relies on simulation outputs, traditional optimization via simulation methods can be applied. We will detail those methods in Section~\ref{sec:SimOpt}. However, many optimization via simulation methods do not consider a sequential decision making setup, i.e., multi-stage, because of the complexity of the convoluted value function structure without a closed form. Therefore, previous literature reviews mainly focus on simulation techniques and convergence proofs in a single-stage or a two-stage setting and do not utilize the potential structure properties of the value functions within a sequential decision process. We believe it is also helpful to discuss the multi-stage optimization under uncertainty solution techniques so that we can better combine them with the simulation methodology and obtain new algorithms. We start our review on the multi-stage optimization under uncertainty methods with a multi-stage extension of model~\eqref{prob:stoch_prob} as follows:
    \begin{align} \label{prob:stoch_prob_multi_single}
        \min_{x_1 \in \cX_1} \; c_1(x_1) + \bE_{\xi_2}[\min_{x_2 \in \cX_2(x_1, \xi_2)} \; c_2(x_2) + \bE_{\xi_3}[\cdots \bE_{\xi_T} [\min_{x_T \in \cX_T(x_1,\dots,x_{T-1}, \xi_2,\dots, \xi_T)} c_T(x_T)]]].
    \end{align}
    We can also present model~\eqref{prob:stoch_prob_multi_single} in a recursive way:
    \begin{align} \label{prob:stoch_prob_multi_recursive}
        \min_{x_1 \in \cX_1} \; c_1(x_1) + \bE_{\xi_2}[f_2(x_1,\xi_2)],
    \end{align}
    where the value function \(f_t(x_1,\dots,x_{t-1};\xi_2,\dots,\xi_t)\) is an optimization problem over the future time horizon \(\{t, \dots, T\}\):
    \begin{align} \label{prob:stoch_prob_multi_vf}
        f_t(x_1,\dots,x_{t-1};\xi_2,\dots,\xi_t) = \min_{x_t \in \cX_t(x_1,\dots,x_{t-1},\xi_2,\dots,\xi_t)} c_t(x_t) + \bE_{\xi_{t+1}}[f_{t+1}(x_1,\dots,x_{t};\xi_2,\dots,\xi_{t+1})]
    \end{align}
    
    At stage \(t = 1\), the first-stage decision \(x_1\) is made based on the known constraint information \(\cX_1\) and cost function \(c_1\). However, the decision \(x_1\) aims to optimize more than just the first-stage cost but the expected future cost. At each following stage \(t = 2,\dots, T\), the uncertainty \(\xi_{t}\) realizes, which helps define the feasible region of the stage-\(t\)'s decision and cost function. The multi-stage stochastic program (or its variants of risk-averse or distributionally robust optimization) often relies on a known and well-structured closed form of \(f\), for example, \(f(x,\xi)\) is a convex linear program in a stochastic programming setting, optimization methods focus on how to approximate the function \(f\). We discuss recent advances in optimization under uncertainty, especially under the multi-stage setting, in Section~\ref{sec:StochOpt}.

    When we want to extend the sequential decision making process to consider black-box future value functions \(f\), there are some main conflicts between the current optimization via simulation and sequential decision making. The sequential decision making process often relies on a rigid framework and assumptions, which cannot approximate any process to an arbitrary precision. The optimization via simulation methods focus more on the statistical analysis and do not fully utilize the structures of the developed optimization algorithms, e.g., Benders decomposition and dual decomposition algorithms.
    
    Therefore, in this paper, we review the recent progress in optimization via simulation and sequential decision processes in the hope that we can develop new methods to combine their advantages and solve complicated multi-stage problems. We will further discuss potential opportunities to combine the two research streams via Bayesian optimization and other meta-models in Section~\ref{sec:SeqOpt}.

\section{Review of Traditional Optimization via Simulation Work} \label{sec:SimOpt}
    Based on the structure and property of the feasible region \(\cX\), we may divide optimization via simulation (OvS) problems into three categories.
    
    \begin{itemize}
        \item If \(\cX\) has a small number of solutions (often less than 100) and the decision vector \(x\) may be numerical or categorical, then we may simulate all solutions and select the best among them. This problem is known as the ranking-and-selection (R\&S) problem.
        \item If \(\cX\) is a convex subset of \(\mathbb{R}^d\) and \(x\) is a vector of continuous decision variables, this problem is known as the continuous OvS (COvS) problem.
        \item If \(\cX\) is a set of \(d\)-dimensional integers and \(x\) is a vector of discrete-valued variables, this problem is known as the discrete OvS (DOvS) problem.
    \end{itemize}

\subsection{Ranking \& Selection}
    There are many reviews available on R\&S, such as~\cite{bechhofer1995design, goldsman1998comparing, kim2006selecting}. Solution procedures can be categorized as two-stage procedures~\cite{dudewicz1975allocation,rinott1978two}, two-stage procedures with screening~\cite{nelson2001simple}, and fully-sequential procedures~\cite{kim2001fully}, all of which can guarantee a desired probability of correct selection. Most procedures above employ a concept of indifference zone (IZ), where a difference is considered significant if larger than a pre-specified indifference-zone parameter. The probability of correct selection guarantee in the IZ approach corresponds to the probability of selecting the true best, subject to the condition that the mean of the true best is better than the mean of all other alternatives by at least the indifference-zone parameter. Thus, this is based on a worst-case performance metric. This worst-case approach can provide frequentist guarantees for correct selection but might require more samples to obtain that guarantee than may be practically implementable. As such, IZ R\&S procedures can be statistically conservative.
    
    There are also R\&S procedures based on average-case performance metrics that sample in a highly sequential manner. The goal is either to maximize evidence for correct selection subject to a constraint on the sampling budget or to reach a level of evidence for correct selection with the fewest expected number of samples. Several typical approaches use an average-case analysis rather than a worst-case bound of the IZ approach, such as optimal computing budget allocation (OCBA) proposed by~\cite{chen2000simulation} and expected value of information (EVI) proposed by~\cite{chick2001new}.
    
    As mentioned above, many existing R\&S procedures are designed to solve small- or medium-scale problems. Researchers begin investigating parallel computing for R\&S problems, such as \cite{luo2015fully,hunter2017parallel,ni2017efficient,kaminski2018parallel,zhong2022knockout,zhong2022speeding}. Meanwhile, some research problems that expand classical R\&S from different perspectives have emerged, e.g., considering
    multiple performance measures by multi-objective R\&S~\cite{batur2018methods, feldman2018score,applegate2020multi}, taking the input uncertainty into account~\cite{gao2017robust,wu2019fixed,fan2020distributionally}, and treating the performance measure as a function of the underlying contexts, which is called R\&S with covariates or contextual R\&S~\cite{shen2017ranking,gao2019selecting,shen2021ranking,ding2022knowledge}. 

\subsection{Gradient Simulation}
    For optimization problems with continuous-valued decision variables, the availability of gradient information can significantly improve the effectiveness of solution algorithms. However, in the stochastic setting, since the outputs are themselves random, deriving stochastic gradient estimators can be challenging.

    The stochastic approximation (SA) algorithm is often used to solve $\mathrm{COvS}$ problems iteratively. Each iteration takes the following form of update:
    $$
    x_{n+1}=\Pi_{\cX}\left[x_n-a_n \widehat{\nabla} f\left(x_n\right)\right],
    $$
    where $x_n$ is the solution found at iteration $n, \widehat{\nabla} f\left(x_n\right)$ is an estimate of the gradient $\nabla f\left(x_n\right),\left\{a_n\right\}$ is a sequence of positive real numbers called the gain sequence, and $\Pi_{\cX}$ denotes a projection operator on to the set $\cX$. We can view SA as the stochastic counterpart to the steepest descent in deterministic optimization. SA is first introduced in~\cite{robbins1951stochastic} to solve noisy root-finding problems and is later applied to the setting of stochastic optimization by solving for the zero of the gradient. Kiefer-Wolfowitz (KW) SA algorithm by~\cite{kiefer1952stochastic} addresses the gradient-free setting.

    We can divide the approaches to stochastic gradient estimation into two main categories, which are indirect and direct.
    
    When the simulation model is treated as a black box, $\nabla f(x)$ may only be estimated using finite-difference (FD) approximations. Let $e_i$ denote the $i$ th column of a $d \times d$ identity matrix. Then a forward FD estimator of $\nabla f(x)$ is
    $$
    \widehat{\nabla} f(x)=\frac{1}{h}\left[\left(\begin{array}{c}
    Y\left(x+h e_1\right) \\
    Y\left(x+h e_2\right) \\
    \vdots \\
    Y\left(x+h e_d\right)
    \end{array}\right)-Y(x)\right],
    $$
    where $h$ is set as a small positive value. The forward FD estimator is biased, and the bias is of $O(h)$. When $Y\left(x+h e_i\right)$ and $Y(x)$ are independent, the variance is of $O\left(h^{-2}\right)$. When the dimension of $x$ is large, the KW SA algorithm is inefficient. To overcome the difficulty, Ref.~\cite{spall1992multivariate} proposes a simultaneous perturbation (SP) FD estimator called SPSA that uses only two simulation runs to obtain an estimate of $\nabla f(x)$ regardless of the dimension of $x$. 

    When the internal structure of the simulation model is known, we may be able to use this knowledge to design more efficient gradient estimators. Among many different methods, the most commonly used are perturbation analysis (PA) and the likelihood ratio/score function (LR/SF) method. Ref.~\cite{ho1983perturbation} proposes PA based on the fact that
    $$
    \nabla f(x)=\nabla \mathrm{E}[Y(x)]=\mathrm{E}[\nabla Y(x)],
    $$
    when the random function $Y(x)$ is stochastically Lipschitz continuous and differentiable with probability one. Then, if $\nabla Y(x)$ can be observed in the simulation, PA estimates $\nabla f(x)$ by $\nabla Y(x)$. Several methods have been proposed for problems where $Y(x)$ is not stochastically Lipschitz continuous, including the smoothed perturbation analysis in~\cite{fu1997conditional} and \cite{gong1987smoothed}, and the pathwise method in \cite{hong2010pathwise}.
    
    Ref.~\cite{glynn1990likelihood} and \cite{reiman1989sensitivity} first introduce the LR/SF method. We can often represent $Y(x)$ by $Y(x)=h(Z)$, where $Z$ is a vector of random variables generated in the simulation and $h(Z)$ is the performance measure in which we are interested. Let $f_z(z, x)$ denote the probability density of $Z$. Then, under some regularity conditions,
    $$
    \begin{aligned}
    \nabla f(x) & =\nabla \mathrm{E}[h(Z)]=\nabla \int h(z) f_z(z, x) d z=\int h(z) \nabla_{x} f_z(z, x) d z \\
    & =\int h(z) \nabla_{x} \log \left[f_z(z, x)\right] f_z(z, x) d z=\mathrm{E}\left\{h(Z) \nabla_{x} \log \left[f_z(Z, x)\right]\right\}.
    \end{aligned}
    $$
    If $f_z(z, x)$ is known, the LR/SF method estimates $\nabla f(x)$ by $h(Z) \nabla_{x} \log \left[f_z(Z, x)\right]$. Similar to the PA estimator, the LR/SF estimator is also unbiased and only requires a single simulation run to compute. Also, it does not require $Y(x)$ to be stochastic Lipschitz continuous. Therefore, its applicability is wider than the PA estimator. Despite the benefits of LR/SF, it has poor variance properties, its variance can be shown to scale linearly with the scale of $Z$ \cite{mcgill2012efficient}.
    
    In recent years, there have been advancements in gradient estimation techniques. Using FD schemes, it takes $O(d)$ function evaluations to obtain just one gradient evaluation, which can be expensive for large $d$ even if parallel computations can be performed. Thus, randomized finite difference methods (RFD) are proposed to estimate gradients using an arbitrary number of function value samples and taking an average of directional derivative estimates along random directions, e.g.,~\cite{fazel2018global,berahas2019derivative,akhavan2022gradient,lobanov2023gradientfree}. In addition, several variance reduction techniques for estimators are proposed, such as control variates \cite{tucker2017rebar,baker2019control,lievin2020optimal} and alternative variational objectives \cite{mnih2016variational,masrani2019thermodynamic}. 

\subsection{Global Methods}
    For discrete optimization via simulation, we use a set of general problem-solving strategies named meta-heuristics. Meta-heuristic algorithms attempt to locate decent feasible solutions to an optimization problem efficiently. To this end, they evaluate potential solutions and perform a series of operations on them to find different, better solutions. We would like to point out that there are extensive works conducted on global methods in the domains of machine learning, but we exclusively focus on presenting the commonly employed global methods in the simulation field due to the scope of this review.
    
    Random search methods involve sampling a set of feasible system designs in each iteration, conducting simulations of the sampled designs to estimate the performance of these designs, and then using the simulation results to decide on what designs should be sampled in the next iteration and on the current estimate of the optimal system design. Consequently, the class of optimization algorithms under consideration is broad enough to be applied to solve both deterministic and stochastic optimization problems with either discrete or continuous decision parameters (or both). There are many random search algorithms for DOvS problems, including the stochastic ruler method of \cite{yan1992stochastic}, the random search methods of \cite{andradottir1995method}, the simulation annealing algorithm of \cite{alrefaei1999simulated}, the stochastic comparison method of \cite{gong2000stochastic}, the nested partitions algorithms of \cite{pichitlamken2003combined} and \cite{shi2000nested}, the model reference adaptive search (MRAS) of \cite{hu2007model}, the COMPASS algorithms of \cite{hong2006discrete} and \cite{hong2007framework}. One desirable feature of random search methods for discrete simulation optimization is that such methods can be shown to converge almost surely to the set of global optimal solutions of the underlying (discrete) simulation optimization problem under very general conditions, see~\cite{andradottir1996global,andradottir2006simulation}. In addition, random search methods have other desirable features, including the ability to move rapidly within the feasible region to identify areas worthy of further investigation, to focus the search on desirable subsets of the feasible region that have not been extensively explored, and to adapt the information available in the current optimization problem; see~\cite{andradottir2006overview} for a more detailed review for such features.
    
    Another strategy to find good solutions is to use information on the past search progress and record this information in memory structures. Metaheuristics that use this strategy are commonly grouped under the umbrella term tabu search algorithms~\cite{glover1989tabu,glover1990artificial,glover1996critical}. A tabu list records the last encountered solutions and forbids these solutions from being revisited as long as they are on the list. Alternatively, the tabu list may also record the last moves that have been made to prevent them from being reversed. The tabu list (tenure) is constructed specifically for the problem structure, and sometimes, it might be hard to design such a list to achieve efficiency.
    
    Simulated annealing (SA) is an iterative search method inspired by metal annealing \cite{kirkpatrick1983optimization}. The algorithm starts with an initial solution and performs a stochastic partial search of the state space using appropriate perturbations and evaluation functions. In each iteration $n$, a random candidate solution $x^{\prime}_n$ is selected from the neighborhood of the current solution $x_{n-1}$. Uphill moves are accepted with a probability $\exp{(-\left[f(x^{\prime}_n)-f(x_{n-1})\right] / T)}$ where $f(\cdot)$ is the objective function value (to be maximized) of the solution between brackets and $T$ is an endogenous parameter called temperature. The probability of acceptance of uphill moves decreases as $T$ decreases. At high temperatures, the search is almost random, while at low temperatures, the search becomes almost greedy, i.e., only good moves are accepted \cite{kirkpatrick1983optimization}. In addition to requiring that the temperature $T$ be strictly positive for all $n$, most of the simulated annealing literature assumes that $T_n\to 0$ as $n \to\infty$ at a logarithmic rate, e.g.~\cite{hajek1988cooling}.

    Genetic algorithms \cite{holland1992genetic}, particle swarm optimization \cite{kennedy1995particle}, and ant colony algorithms \cite{dorigo1992optimization} are also widely employed global methods in simulation. These methods utilize distinct strategies to explore the search space and find optimal solutions. Genetic algorithms simulate the process of biological evolution, employing operations such as selection, crossover, and mutation to explore the search space and find optimal solutions. Particle swarm optimization algorithms simulate collective behavior observed in bird flocks or fish schools, continuously adjusting particle positions and velocities to search for the optimal solution. Ant colony algorithms simulate the foraging behavior of ants, mimicking ant movement and the release of pheromones in the search space to search for the optimal solution. The selection of a suitable method depends on the specific characteristics and requirements of the problem at hand. 
    
\section{Review of Traditional Sequential Decision Process Work} \label{sec:StochOpt}
    The origin of optimization under uncertainty models can be traced back to Dantzig's 1955 paper to incorporate random variables as model parameters in linear programs~\cite{dantzig1955linear}. Since then, many variants have been developed, such as stochastic programming~\cite{birge2011introduction, shapiro2021lectures}, robust optimization~\cite{ben2009robust}, risk-averse optimization~\cite{ruszczynski2013advances}, chance-constrained optimization~\cite{charnes1959chance, pagnoncelli2009sample}, and distributionally robust optimization~\cite{rahimian2019distributionally}. Many of the papers focus on a two-stage setting, where the decision-makers optimize the planning decision in the first stage, observe the uncertainty, and then make the second-stage adjustment~\cite{shapiro1998simulation, schultz1996two, zeng2013solving, liu2016decomposition}. The models listed above are generally hard to directly solve via commercial modeling languages and solvers and thus require transformation, approximation, and decomposition to solve them. It is impossible to list all the methods to solve optimization under uncertainty problems, so we list four fundamental ones here: 1) transformation: for example, robust counterparts derivation for robust optimization models with row-wise uncertainty~\cite{bertsimas2004price, bertsimas2011theory}; integer programming~\cite{luedtke2008sample} and nonlinear programming~\cite{pena2020solving} reformulations for the chance-constrained programs; 2) value function approximation via inner and outer approximation methods~\cite{philpott2013solving, shapiro1998simulation, nemirovski2009robust}; bounding techniques based on relaxations~\cite{huang1977bounds, gassmann1986tight}; 3) decomposition methods such as L-shaped (Benders) cutting-plane methods~\cite{van1969shaped, geoffrion1972generalized} and dual decomposition methods~\cite{caroe1999dual, rockafellar1991scenarios}; 4) stochastic approximation that uses stochastic gradient information for the descent direction~\cite{lan2012validation}, where we essentially consider the stochastic program as a COvS problem but we do have structural information about the second-stage value function to obtain the gradient easily.

    When we extend the optimization under uncertainty models to a multi-stage setting, the solution methods become significantly complicated, and some previously described two-stage transformations no longer apply. For example, generating a stochastic gradient for a particular stage's value function is challenging because the following stages' value functions have not been exactly characterized. Only recently Ref.~\cite{lan2021dynamic} develops a stochastic approximation algorithm that can generalize for stochastic programs with more than three stages, but it requires a nested process of sub-gradient evaluation, which does not scale as the number of stages grows. Therefore, it is generally hard for the iterative algorithms widely used in nonlinear optimization and COvS to satisfy the objective of solving multi-stage problems, and we put our focus on the recent advances in models and solution methods to characterize the form of value functions. For multi-stage stochastic programming problems, most solution methods focus on the value function approximation since the expression \(\bE_{\xi_{t+1}}[f_{t+1}(x_1,\dots,x_{t};\xi_2,\dots,\xi_{t+1})]\) does not provide a directly solvable formulation. Ref.~\cite{pereira1991multi} first proposes the stochastic dual dynamic program (SDDP), which iteratively generates a valid lower approximation at the sampled solutions. For a multi-stage stochastic program, the expected value functions are approximated from below by linear cutting-planes and we are able to solve the following problem: 
    \begin{subequations}
        \begin{align} \label{prob:stoch_prob_multi_vf_lower}
            \underline{f}_t(x_1,\dots,x_{t-1};\xi_2,\dots,\xi_t) = \min_{x_t \in \cX_t(x_1,\dots,x_{t-1},\xi_2,\dots,\xi_t)} \quad & c_t(x_t) + \theta_{t+1} \\
            \st \quad & \theta_{t+1} \geq \pi_{t+1,l}^\top x_t + v_{t+1, l} \qquad l = 1,2,\dots \label{cons:linear_cuts}
        \end{align}
    \end{subequations}
    
    The cut coefficients, \(\pi_{t+1,\cdot}\) and \(v_{t+1,\cdot}\), are obtained from solving lower level optimization problems \(\underline{f}_{t+1}\), i.e., \emph{backward pass}, at the solutions to the sample path problem \(\hat{x}_1,\dots, \hat{x}_T\), simulated and obtained in \emph{forward pass}. If we assume convexity in the value function of every stage, SDDP guarantees asymptotic convergence~\cite{philpott2008convergence, girardeau2014convergence}, and Ref.~\cite{lan2022complexity} provides a thorough discussion about the complexity of SDDP. Because of the computational viability of SDDP and its variants, such cutting-plane algorithms have found wide applications in energy system planning~\cite{gjelsvik2010long}. In recent years, researchers have incorporated more risk-averse aspects in the multi-stage modeling, including the multi-stage distributionally robust optimization~\cite{philpott2018distributionally, duque2020distributionally, arora2022data} and the multi-stage risk-averse optimization~\cite{philpott2013solving, shapiro2009time}. Under a risk-averse setting, the time consistency issue has been central to guarantee an implementable policy~\cite{shapiro2009time, homem2016risk, ruszczynski2010risk}, which requires careful modeling to achieve a composite one-step risk measure structure. For the multi-stage robust optimization, previous works mainly focus on adaptive decision rules or finite decision space, as the multi-stage robust optimization model is convex but infinite-dimensional~\cite{ben2004adjustable, georghiou2015generalized}. Ref.~\cite{georghiou2019robust} develops a new cutting-plane algorithm similar to SDDP but maintains both upper and lower bounds for the value functions. Another important direction is the attempt to address the nonconvexity of the value function. If the value function \(f_t\) is no longer convex, the linear cut~\eqref{cons:linear_cuts} may not provide a tight approximation. For the multi-stage stochastic mixed-integer programs, Ref.~\cite{zou2016nested} proposes Lagrangian cuts, which can be tight and valid if the state variables \(x_t\) are binary.

\begin{figure}
        \centering
        \includegraphics{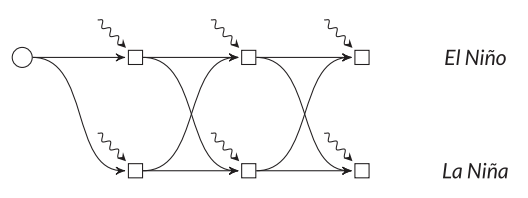}
        \caption{An illustration of the Markovian policy graph in~\cite{dowson2020policy}.}
        \label{fig:policy_graph}
    \end{figure}

    Another important stream of models dealing with multi-stage structured uncertainty is the Markov decision process. Solution methods consist of exact dynamic programming techniques, such as value iteration and policy iteration, and approximate dynamic programming techniques where the value functions are approximated by a combination of basis functions~\cite{keller2006automatic}, heuristic functions~\cite{bertsimas2002approximate}, or \(Q\) functions characterized by machine learning models~\cite{busoniu2017reinforcement}; see~\cite{powell2007approximate} for a detailed description of dynamic programming and approximate dynamic programming and~\cite{sutton2018reinforcement} for Q-learning and other reinforcement learning techniques. Recent developments, first proposed in \cite{philpott2012dynamic}, started to integrate the Markov chain probabilistic structure with a value function iterative approximation algorithm, such as the policy graph definition of the uncertainty~\cite{dowson2020policy}, shown in Figure~\ref{fig:policy_graph}, and the partially observable state in~\cite{dowson2020partially}.

\section{Towards Sequential Decision Process via Simulation} \label{sec:SeqOpt}
    In many real-world industrial applications, the process is too complicated to model with high fidelity, and it is usually costly to construct a full-scale simulation model. In addition, we may need to run such simulation models multiple times for each system design, i.e., candidate solution, in an optimization model to evaluate the expected function value. Even with a successful estimation of the function value (the zeroth order information), we are not able to directly use it in an optimization algorithm similar to those listed in Section~\ref{sec:StochOpt}. Therefore, we need approximation models that can accurately characterize the simulation process, are fast to calculate, and easily capture the functional relationship between the simulation output and the design points. When optimizing system design decisions, we can quickly deploy approximation models in place of the complex simulation model to guide the design into a close neighborhood of optimality. Such approximation is also called \emph{meta-model} or \emph{surrogate model}, and it plays a significant role in many engineered systems' optimization~\cite{hussain2015surrogate, barton2006metamodel, do2022metamodel}.

    Many heuristics have been developed for a specific engineering problem, but there is very scarce literature about the theory considering each stage's transition mechanism and objectives as black-box functions. The difficulty lies in globally approximating the stagewise value functions and embedding the approximation in the optimization problem such that we can search for the optimal solution utilizing known optimization algorithms. We list three different types of meta-models and discuss their advantages and limitations. Some of the listed meta-models are based on local first- or second-order approximations, and thus, we need to fit a series of approximations in a nested optimization scheme, which may not be efficient if the number of stages is large. Despite this, we believe it is valuable to cover the local meta-model methodologies, the philosophy of which may inspire the development of new algorithms.

    \subsection{Response Surface Methodology (RSM)} \label{subsec:rsm}
    The \emph{response surface methodology (RSM)} is a widely applied iterative method to approximate the input-output simulation process~\cite{kleijnen2008response,myers2016response}. In each iteration, RSM constructs a local first-order or second-order polynomial regression model to approximate the simulation output function locally in a neighborhood of the incumbent solution. Suppose we denote the value function to approximate as \(f(x)\), e.g., the expected future value function in a multi-stage stochastic program. Without loss of generality, we can write a second-order model to approximate \(f(x)\) locally:
    \begin{equation}
        g(x) = \beta_0 + \beta_1^\top x + x^\top B_2 x,
    \end{equation}
    where \(\beta_0, \beta_1\) and \(B_2\) represent coefficients for the zeroth-, first-, and second-order terms (we could set \(B_2 = 0\) to yield a first-order meta-model). Given a local region to explore, we first select input design points, \(\hat{x}_1, \hat{x}_2, \dots\), based on design of experiment techniques and obtain the estimated expected simulation responses via multiple rounds of simulations executed at each solution, \(f(\hat{x}_1), f(\hat{x}_2), \dots\)~\cite{kleijnen2017design}. The coefficients are fitted at the simulation design points by minimizing the least square errors between the estimation \(f(\hat{x})\) and the prediction \(g(\hat{x})\). The local model provides a gradient descent, and we can perform one step of the gradient descent method, with the step size determined by binary search or backtracking~\cite{nocedal1999numerical,safizadeh1994optimization}. At the beginning of the RSM procedure, we usually fit a first-order meta-model to quickly converge to a neighborhood that is likely to contain the optimum. Towards the end of the RSM process, a second-order polynomial regression is fitted to capture the functional property more accurately at the local minimum. We can run multiple parallel searching processes starting at different input points to avoid being stuck at one local optimum. Towards the end of the RSM procedure, Ref.~\cite{kleijnen2014response} discusses the test of KKT conditions at the obtained design point to check its optimality. RSM can also be extended to include multiple random responses, selecting one as the objective and others to formulate constraints~\cite{kleijnen2008response}.

    \subsection{Stochastic Kriging} \label{subsec:kriging}
    Stochastic kriging models are another type of general purpose meta-models~\cite{cressie1990origins}; see~\cite{kleijnen2009kriging} for a detailed review of kriging models used in simulation and prediction. It is initially proposed to address the spatial interpolation problem in the field of mining geology but later used to fit a function with the sampled data. In the fitting process, the linear or low-order polynomial regression models tend to have an issue of fitting well locally but generating a non-robust mapping from the simulation output to the meta-model globally. The higher-order nonlinear models work well for certain applications, such as queuing systems, but require domain knowledge and specialized fitting algorithms and may overfit~\cite{ankenman2008stochastic}. The use of kriging in the stochastic simulation is first introduced in~\cite{van2003kriging}: the key idea of kriging is to assume that we can calculate all responses with the sum of a deterministic mean function \(f(x)^\top \beta\), a random field \(M(x)\) that exhibits spatial correlation, e.g., a Gaussian process, and a noise \(\epsilon(x)\):
    \begin{equation} \label{eqn:kriging}
        Y(x) = f(x)^\top \beta + M(x) + \epsilon(x).
    \end{equation}
    We observe \(n_j\) function values \(Y_{ji}, i = 1,\dots,n_j\) at a set of selected design points \(\{\hat{x}_j\}\) for \(j = 1,\dots, k\). The predictor \(\hat{Y}\) at a new design point \(x_0\) is a linear function of the mean simulated values \(\bar{Y}_j = \sum_{i = 1}^{n_j} Y_{ji}\):
    \begin{equation}
        \hat{Y}(x_0) = \lambda_0(x_0) + \sum_{j = 1}^{k} \lambda_j(x_0) \bar{Y}_j.
    \end{equation}
    The coefficients \(\lambda_0,\dots,\lambda_k\) correspond to the solution that minimizes the mean squared error of the predictor assuming the structure of~\eqref{eqn:kriging}; the fitting process for a Gaussian process \(M(x)\) with a power exponential kernel can be found in~\cite{barton2006metamodel}.
    
    The fitted kriging model can drive the search process for the optimal solution. Ref.~\cite{van2003kriging} describes how we can fit a kriging model to the simulation data to search for an optimal design. It involves iteratively simulating new samples according to the fitted kriging model, re-fitting the kriging model according to the newly acquired simulated data, and determining whether KKT conditions are satisfied.

    \subsection{Gaussian Process and Bayesian Optimization} \label{subsec:bayesian}
    Gaussian processes have been widely used for approximating the nonlinear functions from dynamical systems because of their readiness to serve as the prior probability distribution in iterative Bayesian updates; see an introduction for Gaussian process regression in~\cite{beckers2021introduction}. A Gaussian process \(M(x)\) can be considered a distribution over functions with a continuous domain \(x \in \cX\), defined by a mean function \(\bE[M(x)] = \mu(x)\) and a covariance function \(\Sigma(x_1,x_2) = cov(M(x_1),M(x_2))\) to characterize the covariance between two random variables obtained from the Gaussian process at any two points \(x_1\) and \(x_2\) in the domain. Gaussian processes is a global meta-model, which yields a continuously differentiable function for prediction. The smoothness property of Gaussian processes can effectively reduce the mean square error compared to other polynomial regression models in trajectory modeling~\cite{schulz2018tutorial}.

    Since the conjugate of a normal distribution is still normal, we can update the Gaussian process estimation by the Bayes law. Suppose we have simulated at a series of design points \(x_1, \dots, x_n\) and obtained their responses \(f(x_1), \dots, f(x_n)\). By proposing a mean function \(\mu_0(x)\) and a covariance function (or kernel) \(\Sigma_0(x,y)\), which should grant a large positive correlation to two points \(x\) and \(y\) that are close in the input space, we yield a prior distribution on \([f(x_1), \dots, f(x_n)]\) as:
    \begin{equation}
        f(x_{1:n}) \sim Normal(\mu_0(x_{1:n}), \Sigma_0(x_{1:n}, x_{1:n})),
    \end{equation}
    where we follow the notation in~\cite{frazier2018tutorial} to let the subscripts \(1:n\) represent the vectorized values. Suppose in the \(n+1\) observation, we would like to acquire the information at a new design point \(x\). We can calculate the conditional distribution of \(f(x)\) based on the observed data:
    \begin{subequations}
        \begin{align}
            f(x)|f(x_{1:n}) & \sim Normal(\mu_n(x), \sigma_n^2(x)) \\
            \mu_n(x) & = \Sigma_0(x, x_{1:n})) {\Sigma_0(x_{1:n}, x_{1:n}))}^{-1} (f(x_{1:n} - \mu_0(x_{1:n})) + \mu_0(x) \\
            \sigma_n^2(x) & = \Sigma_0(x,x) - \Sigma_0(x,x_{1:n}) \Sigma_0(x_{1:n},x_{1:n})^{-1} \Sigma_0(x_{1:n},x).
        \end{align}
    \end{subequations}
    \cite{beckers2021introduction} list a detailed derivation process of this conditional distribution, i.e., the posterior distribution \(f(x)|f(x_{1:n})\) and discussed methods to select the appropriate mean function and kernel function, e.g., Mat\'{e}rn kernel and squared exponential kernel. There are many approaches to optimize hyper-parameters in the kernel function, such as maximum a posteriori (MAP) estimate~\cite{gelman1995bayesian}, log marginal likelihood approach~\cite{williams2006gaussian}, and fully Bayesian approach~\cite{rossi2021sparse}.~\cite{williams2006gaussian} also discuss more efficient algorithms to perform this Bayesian update step.

    Based on the Bayesian update of the Gaussian process, Ref.~\cite{kushner1964new} proposes a search algorithm, i.e., \emph{Bayesian optimization}, which is later popularized by the foundational work of Efficient Global Optimization (EGO) in~\cite{jones1998efficient}. Ref.~\cite{brochu2010tutorial},~\cite{shahriari2015taking},~\cite{frazier2018tutorial}, and many others thoroughly review the Bayesian optimization approach, stating the situation for which Bayesian optimization is suitable, including low dimension for the design point \(x\), a feasible region with a simple structure, and continuity. Ref.~\cite{frazier2018tutorial} summarizes the Bayesian optimization algorithm by the following steps:
    \begin{enumerate}[label=Step~\arabic*]
        \itemsep0.3em
        \item Initialize with a Gaussian process prior on \(f\), including the definition of the mean and kernel functions \(\mu_0\) and \(\Sigma_0\);
        \item Observe \(f\) at \(n_0\) design points and set \(n = n_0\);
        \item Update the posterior probability distribution on \(f\) based on the observation;
        \item Select the maximizer of the acquisition function, \(x_{n+1}\), and observe \(f(x_{n+1})\);
        \item Increment \(n\) and repeat Step 3-4 until the iteration limit.
    \end{enumerate}
    
    In each step, the random variable that represents the future value function \(f(x,\xi)\) is approximated by a Gaussian process posterior. The posterior Gaussian process is then used as the prior input for the next iteration. The acquisition function guides the search algorithm to the next candidate design point \(x_{n+1}\), where we obtain new samples to update the estimation. The acquisition function is chosen to balance the exploration and the exploitation, a concept commonly practiced in multi-armed bandit problems and reinforcement learning~\cite{gupta2006interplay}. If maximizing the acquisition function leads to a new design point that is too close to the best design point among the evaluated ones, the search algorithm may not be able to escape the local minimum and explore other potential regions. Common acquisition functions include the expected improvement~\cite{jones1998efficient}, knowledge-gradient~\cite{frazier2009knowledge}, entropy search~\cite{hennig2012entropy}, and predictive entropy search~\cite{hernandez2014predictive}. While the expected improvement and its variants are the most widely used, other alternatives provide validity under non-standard assumptions.

    It requires adaptation of vanilla Bayesian optimization algorithms to run in parallel and handle noisy measurements and complicated constraints, which are often necessary for real-world applications and highly related to optimization under uncertainty. For constraints without noise,~\cite{schonlau1998global} define the expected improvement to be zero if a candidate design point is infeasible. Such expected constrained improvement is also independently developed by~\cite{gardner2014bayesian}. Ref.~\cite{letham2017constrained} develops a Bayesian expected improvement step for noisy observations and constraints. The integration step over the posterior of the acquisition function is achieved via a quasi-Monte Carlo approximation. Ref.~\cite{eriksson2021scalable} proposes a state-of-the-art scalable Bayesian optimization method that employs trust regions to confine the search points and transformations for the objective function and constraints for easier feasibility/optimality detection.

    \subsection{Discussion: Bridging between Simulation Meta-models and Sequential Decision Making}
    Most meta-models aim to solve a single-stage optimization problem with a similar form of model~\eqref{prob:stoch_prob}, as it is difficult enough to develop a converging algorithm that simultaneously performs functional estimation and optimum search while providing a statistical guarantee. Few previous works extend the meta-models listed in Section~\ref{subsec:rsm}-\ref{subsec:bayesian} to a sequential, or even simply two-stage, setting. 
    Ref.~\cite{bailey1999response} applies RSM to approximate the objective function of a two-stage stochastic linear program with recourse, with candidate solutions selected by central composite design (CCD) and the expected value function estimated by sample average approximation and Latin hypercube sampling. However, such a method does not aim to solve the two-stage stochastic program but uses RSM to perform a sensitivity analysis on the optimal solution. Recently, Ref.~\cite{zhuo2022rsm} represents the value functions in an approximate dynamic programming setting via RSM and generalized polynomial chaos models for a stochastic energy management problem. Ref.~\cite{xie2019global} develops a two-stage optimization via simulation method. In the recourse problem, the objective function is replaced by a local Gaussian process. The expectation of recourse functions is then approximated by a global Gaussian process. Iterating between the global-local meta-models, the algorithm can guarantee the convergence to the optimum for two-stage discrete stochastic programs.
    
    There are still some challenges to extending the current optimization via simulation methods to a multi-stage setting. Here, we try to list a few points we believe could require more future research:
    \begin{itemize}
        \item Many of the simulation optimization methods can select multiple design points and compare them, then finalize the best decision. We disregard the cumulative objective values at the trial design points as long as the best decision can be obtained. On the other hand, we only have historical data for data-driven sequential decision making and need to optimize the current trial. We believe it is necessary to identify the correct goal before we propose new optimization methods. Do we have a fixed budget for trial-and-error and only care about the optimality of the final decision after the sampling process? Or do we want to find the optimal solution only based on the historical data? Or do we rather want to optimize the long-run performance where each trial's optimal value also counts, but we would also like to reveal more information? Different applications may find one of those three to be the most suitable.
        \item The noise in the meta-models' observations may have a significantly smaller magnitude than the uncertainty in sequential decision making. This may require many more samples to obtain an accurate functional evaluation and makes it hard to achieve any reasonable result when the sample budget or the historical data is limited.
        \item It is hard for meta-models to deal with the nonconvexity and feasibility requirements. The value function is evaluated at every stage by solving an optimization problem with the future value functions embedded in it. It is computationally difficult to formulate the future value functions via Gaussian processes and other meta-models because they most likely will induce nonconvexity and cannot characterize the feasible region of the current stage in which we guarantee to find a feasible solution in all future stages for any uncertainty realization. 
        \item The recursive relationships between stages make the accuracy of meta-model approximation decay fast. One potential idea to improve the accuracy is to use machine learning models as surrogates since they can approximate functions with high precision, especially deep neural networks, but it is still computationally challenging since they may require a mixed-integer program or other nonconvex representation, e.g., in Neur2SP~\cite{dumouchelle2022neur2sp}. This may also face issues when the data points are insufficient.
    \end{itemize}

    \section{Conclusions}
    In this survey, we review the development of methodology in optimization via simulation and sequential decision processes. As pointed out by many previous literature, the fast development of computer simulation models makes them more accurate and applicable for a wide variety of applications, but embedding such high-fidelity simulation processes in sequential decision making still needs more effort, especially in handling the nonconvexity and infeasibility. It is impossible to include every paper with a related topic. Still, we hope this survey can inspire new research to bridge the gap for sequential decision making via simulation.
    
    \bibliographystyle{abbrv}
    \bibliography{Multi-stage_SimOpt}

\end{document}